\documentclass[aos,MSNbibl,nameyear,dvips]{arximspdf}
\usepackage{graphicx}
\doi{10.1214/13-AOS1130} 
\volume{41}
\issue{4}
\pubyear{2013}
\firstpage{1944}
\lastpage{1969}

\makeatletter

\newcommand{\rright}{\right}
\newcommand{\lleft}{\left}
\newtheorem{thmm}{Theorem}[section]
\newtheorem{lem}[thmm]{Lemma}
\newtheorem{prop}[thmm]{Proposition}

\newproclaim{defn}[thmm]{Definition}
\newproclaim{defprop}[thmm]{Definition/Proposition}
\newproclaim{const}[thmm]{Construction}

\newproclaim{exmp}[thmm]{Example}
\newproclaim{assp}[thmm]{Assumption}
\newproclaim{ques}{Question}

\newproclaim{rem}[thmm]{Remark}
\newproclaim{note}{Note}
\newproclaim{case}{Case}

\newcommand{\RR}{\mathbb{R}}
\newcommand{\NN}{\mathbb{N}}
\newcommand{\by}{\mathbf{y}}
\newcommand{\GL}{\operatorname{GL} }
\newcommand{\cG}{\mathcal{G}}
\newcommand{\cS}{\mathcal{S}}
\newcommand{\tcG}{\widetilde{\cG}}
\newcommand{\SG}{\cS_\cG}
\newcommand{\SplusG}{\cS_\cG^+}
\newcommand{\bX}{\mathbf{X}}
\newcommand{\down}{\hspace*{-1pt}\downarrow\hspace{-3pt}}
\newcommand{\downb}{ \down}
\newcommand{\bT}{\mathbf{T}}
\renewcommand{\ggg}{\mathfrak{g}}
\newcommand{\bP}{\mathbf{P}}
\newcommand{\cC}{\mathcal{C}}
\newcommand{\Aut}{\mathrm{Aut}}
\makeatother

\begin{document}
\begin{frontmatter}

\title{Groups acting on Gaussian graphical models}
\runtitle{Groups acting on Gaussian graphical models}

\begin{aug}
\author[A]{\fnms{Jan} \snm{Draisma}\ead[label=e1]{j.draisma@tue.nl}\thanksref{t1}},
\author[B]{\fnms{Sonja} \snm{Kuhnt}\ead[label=e2]{kuhnt@statistik.tu-dortmund.de}}
\and
\author[C]{\fnms{Piotr} \snm{Zwiernik}\corref{}\ead[label=e3]{pzwiernik@berkeley.edu}\thanksref{t1}}
\thankstext{t1}{Supported by Draisma's Vidi grant from the Netherlands
Organisation for Scientific Research (NWO).}
\runauthor{J. Draisma, S. Kuhnt and P. Zwiernik}
\affiliation{TU Eindhoven, TU Dortmund University and CWI Amsterdam}
\address[A]{J.Draisma\\
Department of Mathematics\\
\quad and Computer Science\\
TU Eindhoven\\
PO Box 513 \\
5600 MB Eindhoven\\
The Netherlands\\
\printead{e1}}
\address[B]{S. Kuhnt\\
Faculty of Statistics\\
TU Dortmund University\\
44221 Dortmund\\
Germany\\
\printead{e2}}
\address[C]{P. Zwiernik\\
Department of Statistics\\
University of California Berkeley\\
Berkeley, California 94720\\
USA\\
\printead{e3}}
\end{aug}

\received{\smonth{7} \syear{2012}}
\revised{\smonth{5} \syear{2013}}

%
\begin{abstract}
Gaussian graphical models have become a well-recognized tool for the
analysis of conditional independencies within a set of continuous random
variables. From an inferential point of view, it is important to realize
that they are composite exponential transformation families. We reveal
this structure by explicitly describing, for any undirected graph, the
(maximal) matrix group acting on the space of concentration matrices in the
model. The continuous part of this group is captured by a poset naturally
associated to the graph, while automorphisms of the graph account for
the discrete part of the group. We compute the dimension of the space of
orbits of this group on concentration matrices, in terms of the combinatorics
of the graph; and for dimension zero we recover the characterization by Letac
and Massam of models that are transformation families. Furthermore, we
describe the maximal invariant of this group on the sample space, and
we give a sharp lower bound on the sample size
needed for the existence of equivariant estimators of the concentration
matrix. Finally, we address the issue of robustness of these estimators
by computing upper bounds on finite sample breakdown points.
\end{abstract}

%


\end{frontmatter}

\section{Introduction and results} \label{secintro}

Gaussian graphical models are popular tools for
modelling complex associations in the multivariate continuous case. If
the graph with vertex set $[m]:=\{1,\ldots,m\}$ is complete, then the
general linear group $\GL_m(\RR)$, consisting of all invertible $m
\times m$-matrices, acts on the space of concentration matrices in the
model, as well as on the sample space. The maximum likelihood estimator
(MLE) of the concentration matrix is equivariant with respect to this
group action, but many other equivariant estimators have
been proposed, for example, by \citet
{jamesstein,andersonMVST,lopuhaa91,donohophd,stahel1981robust}. For
smaller graphs,
only some proper subgroup of $\GL_m(\RR)$ will act on the set
of compatible concentration matrices. In this paper, we describe that
subgroup explicitly, and pave the way for its use in designing invariant
tests, (robust) equivariant estimators and improved inference procedures.

Having an explicit group acting on a statistical model has numerous
advantages. This was first pointed out by \citet{fisher1934} in the
context of the
location and scale models, which then led to the notion
of a transformation family, that is, a~statistical model on which a
group acts
transitively. Group actions give rise, for example, to the
study of model invariants and distributional aspects of the maximum
likelihood estimator (MLE) or other
equivariant estimators [see \citet{barndorffMLE,nielsentransformation82,eaton1989,fisher1934,reidconditioning,lehmannromano}].
When a group acts on a model
in a nontransitive manner, the model is sometimes called a composite
transformation family [see \citet{nielsentransformation82}]. In this
case, the model can be decomposed into a family of transformation
models each corresponding to a fixed value of some parameter.

To set the stage, let $\cG=([m],E)$\vspace*{1pt} be an undirected graph with set
of vertices $[m]$ and set of edges $E\subseteq
{[m]\choose2}$. Denote by $\cS_{m}$ the set of symmetric matrices
in $\RR^{m\times m}$ and by $\cS^{+}_{m} \subseteq\cS_m$ the cone of
positive definite matrices. Let $\cS_{\cG}\subseteq\cS_{m}$ denote the
linear space of symmetric matrices whose $(i,j)$ off-diagonal entry is
zero if $\{i,j\} \notin E$, and by $\cS^{+}_{\cG}$ the cone of all
positive definite matrices in $\cS_{\cG}$. As a running example in
this \hyperref[secintro]{Introduction},
we take $\cG$ to be the path $P_{3}\dvtx \stackrel{2}{\bullet}
-\stackrel{1}{\bullet}-\stackrel{3}{\bullet}$. So $\cS_{P_3}$
consists of
all symmetric matrices of the form
\[
\lleft[\matrix{ * & * & * \vspace*{2pt}
\cr
* & * & 0 \vspace *{2pt}
\cr
* & 0 & * }
\rright].
\]

Let $X=(X_{i})_{i\in[m]}$ be a random vector with multivariate normal
distribution $\mathcal{N}(0,\Sigma)$. The \emph{Gaussian graphical model}
is the statistical model for $X$ given by
\[
M(\cG):=\bigl\{\mathcal{N}(0,\Sigma) \mid\Sigma^{-1}\in
\cS^{+}_{\cG
}\bigr\};
\]
so $\cS^+_\cG$ is the space of concentration matrices compatible with
the model [see \citet{lauritzen96}].

The group $\GL_{m}(\RR)$ acts on $\RR^m$ by matrix-vector
multiplication, and this induces an action on $\cS_m$ and $\cS^{+}_{m}$
given by $g\cdot K:= g^{-T}K g^{-1}$---indeed, note that this
is the concentration matrix of $g X$ if $K$ is the concentration matrix
of $ X$.

A leading role in this paper is played by the group
\[
G:=\bigl\{g \in\GL_m(\RR) \mid g \cdot\SplusG\subseteq\SplusG\bigr
\}.
\]
This is a closed subgroup of the Lie group $\GL_m(\RR)$ (see
Section~\ref{secG}). For example, if $\cG$ is the complete graph, then
$G$ is all of $\GL_m(\RR)$. For \emph{any} graph $\cG$, the group $G$
contains the invertible diagonal matrices, which correspond to scaling
the components of $X$. Furthermore, $G$ contains elements coming
from graph automorphisms of $\cG$. Specifically, if $\pi\dvtx [m] \to[m]$
is such an automorphism, then the permutation matrix $g$ with ones on
the positions $(i,\pi(i)),\ i \in[m]$ lies in $G$, since its action on
$\cS_m$ stabilizes the zero pattern prescribed by $\cG$. For our running
example $P_3$, the permutation matrix
\[
\lleft[\matrix{ 1 & 0 & 0 \vspace*{2pt}
\cr
0 & 0 & 1 \vspace *{2pt}
\cr
0 & 1 & 0 }
\rright]
\]
lies in $G$.

\subsection{The group $G$}

Our first result is an explicit description of $G$ in terms of~$\cG$,
and requires the pre-order on $[m]$ defined by
\[
i\preccurlyeq j\quad \mbox{if and only if}\quad N(j)\cup\{j\}\subseteq N(i)\cup\{i\},
\]
where $N(i)=\{j \in[m]\dvtx \{i,j\} \in E\}$ denotes the set of neighbors
of $i$ in $\cG$. So in our running example $P_3$ we have $1
\preccurlyeq
2,3$. Consider the closed subset $G^0$ of $\GL_m(\RR)$ defined by
\[
G^0 = \bigl\{g\in\GL_m(\RR) \mid g_{ij}= 0
\mbox{ for all } j\not\preccurlyeq i\bigr\}.
\]
We show in Section~\ref{secG} that this set is a subgroup of
$\GL_m(\RR)$. For $\cG=P_3, $ it consists of all invertible matrices
of the form
\[
\lleft[\matrix{ * & 0 & 0 \vspace*{2pt}
\cr
* & * & 0 \vspace*{2pt}
\cr
* & 0 & *}
\rright].
\]

\begin{thmm} \label{thmmG}
For any undirected graph $\cG=([m],E)$, the group $G$ is generated by
the group
$G^0$ and the permutation matrices corresponding to the automorphism
group of the
graph $\cG$.
\end{thmm}

For $P_3$, this theorem says that $G$ is the group of all matrices of
the form above, together with all matrices of the form
\[
\lleft[\matrix{ * & 0 & 0 \vspace*{2pt}
\cr
* & 0 & * \vspace*{2pt}
\cr
* & * & 0
} \rright].
\]
The two subgroups of $G$ in Theorem~\ref{thmmG} can have a nontrivial
intersection. For instance, when $\cG$ is the complete graph, the
automorphism group of $\cG$ is contained in $G^0$. In Section~\ref{secG},
we state and prove a more refined statement that gets rid of that
intersection.

\subsection{Existence and robustness of equivariant estimators}

Now that we know explicitly which matrix group $G$ acts on our graphical
model $M(\cG)$, we can use this group to develop classical notions of
multivariate statistics in the general context of graphical models. One
of these notions is that of an equivariant estimator [see, e.g., \citet{eaton1989,schervishbook}]. Let $\bX$ denote the $m\times n$
matrix, whose columns correspond to $n$ independent copies of the
vector $X$. Then an equivariant estimator for the concentration
matrix is a map $T\dvtx (\RR^{m})^n = \RR^{m \times n} \to\SplusG$,
that is, a map from the space of $n$-samples $\bX$ to the parameter space
of the model, that satisfies $T(g\mathbf{x})=g T(\mathbf{x})$ for all
realisations $\mathbf{x}$ of $\bX$. The standard example is the
maximum likelihood estimator (MLE). Indeed, the likelihood
of concentration matrix $K$ given an $n$-sample $\mathbf{x}$ equals the
likelihood of $g\cdot K$ given $g \mathbf{x}$, for any $g \in G$, and
this implies that the MLE is $G$-equivariant. Other equivariant
estimators of the concentration matrix for some special graphical
models have been proposed
in \citet{sun2005estimation}.

For decomposable graphs, the MLE exists with probability one if and
only if $n$ is at least the size of the maximal clique of the given
graph. However, in general, whether the MLE exists, with probability
one, for a given sample size
$n$ and a given graph $\cG$ is a subtle matter; see the recent paper by
\citet{uhler2010} and the references therein. By contrast, the
question whether for a given sample size \emph{any} equivariant estimator
exists, turns out to have a remarkably elegant answer for any graph
$\cG$. To state it,
define the \emph{down set} $\down i$ of an element $i \in[m]$ to be
the set of all $j \in[m]$ with $j \preccurlyeq i$.

\begin{thmm}\label{thGeqexist}
Let $\cG=([m],E)$ be an undirected graph. There exists a $G$-equivariant
estimator $T\dvtx \RR^{m\times n}\rightarrow\cS^{+}_{\cG}$ if and only
if $n\geq
\max_{i\in[m]}|\down i|$.
\end{thmm}

To be precise, when $n$ is at least the bound in the theorem, a
$G$-equivariant $T$ exists that is defined outside some measure-zero set
(in fact, an algebraic subvariety of positive codimension), while if $n$
is smaller than that bound, then not even any partially defined
equivariant map
$T$ exists. For our running example $P_3$, we have $\down1=\{1\}$ and
$\down2=\{1,2\}$ and $\down3=\{1,3\}$, so Theorem~\ref{thGeqexist}
says that an equivariant estimator exists with probability one if and only
if the sample size is at least $2$, which in this case coincides with
the condition for existence of the MLE.

Theorem~\ref{thGeqexist} will be proved in
Section~\ref{secequivariance}, where we also establish upper bounds
on the robustness of equivariant estimators, based on general theory
from \citet{daviesgather2005}.

\subsection{The maximal invariant}

Another classical notion related to a group action on a statistical model
is that of \emph{invariants} on the sample space. In our case, these
are maps $\tau$ defined on $\RR^{m \times n}$, possibly outside some
measure-zero set, that are constant on $G$-orbits, that is, that satisfy
$\tau(g\bX)=\tau(\bX)$ for all $g \in G$. An invariant $\tau$ is
called \emph{maximal} if it distinguishes
all $G$-orbits. In formulas this means that for $n$-samples
$\mathbf{x},\by\in\RR^{m \times n}$, outside some set of
measure zero, the equality $\tau(\mathbf{x})=\tau(\by)$ implies
that there exists a $g \in G$ such that $g\mathbf{x}=\by$. Any invariant
map is then a function of $\tau$.

The relevance of maximal invariants in statistics lies
in the fact that they facilitate inference for the maximum likelihood
estimator in the case of transformation families [see \citet{barndorffMLE,reidconditioning}]. In this case the maximal
invariant is an ancillary statistics that one may chose to condition
on. These ideas can be used also in the case of composite
transformation families, where the inference for the index parameter
$\kappa$ is based on the marginal distribution of the maximal
invariant statistics [\citet{barndorffMLE}, Section~5].

Another important application of the maximal invariant is in the
construction of \emph{invariant tests}
[see \citet{eaton1989,lehmannromano}]. Suppose, for
instance, that we want to test the hypothesis that the distribution of
the multivariate Gaussian random vector $X$ lies in $M(\cG)$ against
the alternative that it does not, and suppose that for the $n$-sample
$\bX=\mathbf{x}$ the test would accept the hypothesis. Then,
since $M(\cG)$ is stable under the action of any $g \in G$, it is natural
to require that our test also accepts the hypothesis on observing
$g\mathbf{x}$. Thus, the test itself would have to be $G$-invariant.

Our result on maximal invariants
uses the equivalence relation $\sim$ on $[m]$ defined by $i \sim j$
if and only if both $i \preccurlyeq j$ and $j \preccurlyeq i$, that is,
if and only if $N(i)\cup\{i\}=N(j)\cup\{j\}$. We write $\bar{i}$ for
the equivalence class of $i \in[m]$ and $[m]/\sim$ for the set of all
equivalence classes.

\begin{thmm} \label{thmmmaxinv}
Let $\cG=([m],E)$ be an undirected graph. Suppose that $n\geq
\max_i|\down i|$. Then the map $\tau\dvtx \RR^{m\times n}\to
\prod_{\bar{i}\in[m]/\sim} \RR^{n\times n}$ given by
\[
\mathbf{x}\mapsto \bigl(\mathbf{x}[\downb i]^T \bigl(\mathbf{x}[
\downb i] \mathbf {x}[\downb i]^T\bigr)^{-1} \mathbf{x}[
\downb i] \bigr)_{\bar{i}\in[m]/\sim},
\]
where $\mathbf{x}[\downb i] \in\RR^{|\downb i|\times n}$ is the submatrix
of $\mathbf{x}$ given by all rows indexed by $\down i$, is a maximal
$G^0$-invariant.
\end{thmm}

The lower bound on $n$ in the theorem ensures that the $| \downb
i| \times| \downb i|$-matrices $\mathbf{x}[\downb i] \mathbf
{x}[\downb i]^T$ are
invertible for generic $\mathbf{x}$, and in particular for $\mathbf
{x}$ outside a
set of measure zero. For the complete graph, Theorem~\ref{thmmmaxinv}
reduces to the known statement that $\mathbf{x}\mapsto\mathbf{x}^T
(\mathbf{x}\mathbf{x}^T)^{-1}
\mathbf{x}$ is a maximal invariant, see Example 6.2.3 in \citet
{lehmannromano}, while for our running example $P_3$
it says that the rank-one matrix $\mathbf{x}[1]^T (\mathbf{x}[1]
\mathbf{x}[1]^T)^{-1} \mathbf{x}[1]$
(recording only the direction of the first row of $\mathbf{x}$) and the
rank-two matrices $\mathbf{x}[1,2]^T (\mathbf{x}[1,2] \mathbf
{x}[1,2]^T)^{-1} \mathbf{x}[1,2]$
and $\mathbf{x}[1,3]^T (\mathbf{x}[1,3] \mathbf{x}[1,3]^T)^{-1}
\mathbf{x}[1,3]$ together form a
maximal invariant for $G^0$.

We stress that Theorem~\ref{thmmmaxinv} gives a maximal invariant under
the subgroup $G^0$, rather than under all of $G$. The proof of this
theorem can be found in Section~\ref{secorbequiv}.

\subsection{Orbits of $G$ on $\SplusG$}

Our final results concern the space $\SplusG/G$ of $G$-orbits in
$\SplusG$. When $M(\cG)$ is a transformation family, this space
consists of a single point and hence has dimension zero. Conversely,
it turns out that when the dimension of $\SplusG/G$ is zero, $M(\cG)$
is a transformation family. By work of \citet{letacmassam2007}, it is known
exactly for which graphs this happens. Our result on $\SplusG/G$ is a
combinatorial expression for its dimension. Rather than capturing that
expression in a formula, which we will do in Section~\ref{secorbitS},
we now describe it in terms of a combinatorial procedure.

\begin{figure}

\includegraphics{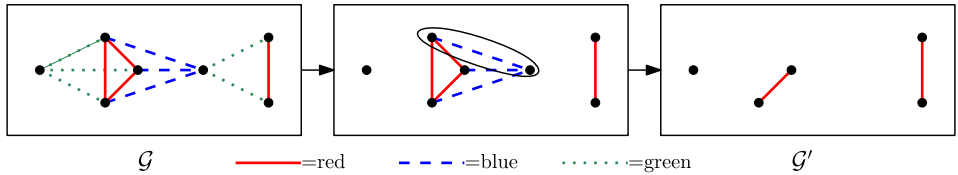}

\caption{An example where the orbit space $\SplusG/G$ has dimension $1$.}
\label{figorbitspacedim}
\end{figure}

Let $\cG=([m],E)$ be an undirected graph. Color an edge $\{i,j\} \in E$
red if $i \sim j$, green if $i \preccurlyeq j$ or $j \preccurlyeq i$
but not both, and blue otherwise. Next delete all green edges from
$\cG$, while retaining their vertices. Then delete the blue edges
sequentially, in each step not only deleting a blue edge but also its
two vertices together with all further blue and red edges incident
to those two vertices. Continue this process until no blue edges are
left. Call the resulting graph $\cG'$; it consists of red edges only. See
Figure~\ref{figorbitspacedim} for an example. One can show that, up to
isomorphism, $\cG'$ is independent of the order in which the blue edges
with incident vertices were removed---though in general it is larger
than the graph obtained by deleting all blue edges, their vertices,
and their incident edges at once.

\begin{thmm} \label{thmmorbitspacedim}
The dimension of $\SplusG/G$ equals the number of blue edges in the
original graph $\cG$ minus the number of red edges in $\cG$ plus the
number of remaining red edges in $\cG'$.
\end{thmm}

In other words, that dimension equals the number of blue edges in $\cG$
minus the number of red edges deleted in the process going from $\cG$
to $\cG'$. This number is nonnegative: indeed, if in some step a blue
edge $\{i,j\}$ is being deleted together with its vertices, then for
each red edge $\{i,k\}$ being deleted along with $i$ there is also
a blue edge $\{k,j\}$ being deleted, and for each red edge $\{j,l\}$
being deleted along with $j$ also a blue edge $\{i,l\}$ is deleted. This
shows, in particular, that $\dim\SplusG/G$ is zero if and only if
$\cG$ has no blue edges, that is, if all edges run between vertices that
are comparable in the pre-order. This is equivalent to the condition
found in \citet{letacmassam2007} for $M(\cG)$ to be a transformation
family; see
Theorem~\ref{thmmmassam} below.

For our running example $P_3, $ the model is a transformation family
and similarly for complete graphs. For an example where $\SplusG/G$
has dimension $1$, see Figure~\ref{figorbitspacedim}. The proof of
Theorem~\ref{thmmorbitspacedim} can be found in Section~\ref{secorbitS}
and in supplementary materials [\citet{supp}].

\subsection*{Organization of the paper}

The remainder of the paper closely follows the structure of this
introduction. First, in Section~\ref{secG} we use
structure theory of real algebraic groups to determine $G$. In
Section~\ref{secequivariance}, we derive necessary and sufficient
conditions for the existence, with probability one, of equivariant
estimators of the concentration (or covariance) matrix, and we give an
upper bound on the robustness of those estimators, measured by the finite
sample breakdown point for generic samples.
In Section~\ref{secorbequiv}, we derive
the maximal invariant of Theorem~\ref{thmmmaxinv}. In
Section~\ref{sectransitive}, we discuss in some detail the case where $G$
acts transitively on $\SplusG$ providing general formula for an
equivariant estimator, after which Section~\ref{secorbitS}
is devoted to our combinatorial formula for the orbit space dimension
in the general case. We conclude the paper with a short discussion. In
the supplementary materials [\citet{supp}], we provide the proof of Theorem~\ref
{thmmorbitspacedim}. We also discuss further results on the
combinatorial structure of the problem that link our work to \citet
{andersson1993}.

\section{The group $G$}\label{seclattice}\label{secG}

Throughout this paper, we fix an undirected graph $\cG=([m],E)$ and
define the group $G$ as in the \hyperref[secintro]{Introduction}:
\[
G:=\bigl\{g \in\GL_m(\RR) \mid g \cdot\SplusG\subseteq\SplusG\bigr
\}.
\]
Note that $G$ is, indeed, a subgroup of $\GL_m(\RR)$: first, if $g,h
\in G$, then $(gh) \cdot\SplusG\subseteq g \cdot\SplusG\subseteq
\SplusG$; and second, if $g \in G$, then since $\SplusG$ linearly
spans $\cS_\cG$, the (linear) action of $g$ must map the linear space
$\cS_\cG$ into itself. Since $g$ is invertible, we then have $g \cdot
\SG= \SG$ (which implies that $g \cdot\SplusG= \SplusG$ holds instead
of the apparently weaker defining inclusion). But then also $g^{-1}
\cdot
\cS_\cG=\cS_\cG$. Finally, the action by $g^{-1}$ preserves positive
definiteness, so that $g^{-1} \cdot\SplusG=\SplusG$, as claimed.

The general linear group $\GL_m(\RR)$ has two natural topologies: the
Euclidean topology, and the weaker Zariski topology in which closed sets
are defined by polynomial equations in the matrix entries. The subgroup
$G$ is closed in both topologies. Indeed, by the above, its elements $g$
are characterized by the condition that $g^T K g \in\SG$ for all $K
\in
\SG$, and this translates into quadratic equations in the entries of $g$.
As a Zariski-closed subgroup of $\GL_m(\RR)$, the group $G$ is a real
algebraic matrix group, and in particular a real Lie group. For basic
structure theory of algebraic groups, we refer to \citet{Borel91}.
In algebraic groups, the Zariski-connected component containing the
identity is always a normal subgroup, the quotient by which is finite. We
first determine the identity component and then the quotient.

\subsection{The identity component}

Observe that the group
\[
\bT^m:= \bigl(\GL_{1}(\RR)\bigr)^{m}\subseteq
\GL_{m}(\RR)
\]
of all invertible diagonal matrices is contained entirely in $G$---indeed,
it just rescales the components of the random vector $X$ and therefore
preserves the original conditional independence statements defining
$M(\cG)$. The group
$\bT^m$ has $2^m$ components in the Euclidean topology, corresponding to
the possible sign patterns of the diagonal entries, but it is connected in
the Zariski topology. For this reason, the Zariski topology is slightly
more convenient to work with, and in what follows our topological
terminology refers to it.

We will use that the connected component of $G$ containing the identity
(the \emph{identity component}, for short) is determined uniquely by
its Lie algebra~$\ggg$. The following lemma helps us determine that Lie
algebra; we use the standard notation $E_{ij}$ for the matrix that has
zeroes everywhere except for a one at position $(i,j)$.

\begin{lem}\label{lemifdiag}
Let $H\subseteq\GL_{m}(\RR)$ be a real algebraic matrix group containing
the group $\bT^m$. Then the Lie algebra of $H$ has a basis consisting of
matrices $E_{ij}$ with $(i,j)$ running through some subset $I$ of $[m]
\times[m]$. Moreover, the set $I$ defines a pre-order on $[m]$ in the
sense that $(i,i)$ lies in $I$ for all $i \in[m]$ and that $(i,j),(j,k)
\in I \Rightarrow(i,k) \in I$. Conversely, the $E_{ij}$ with $(i,j)$
running through any set $I \subseteq[m] \times[m]$ defining a pre-order
on $[m]$ span the Lie algebra of a unique closed connected subgroup
of $\GL_m(\RR)$ containing $\bT^m$, namely, the group of all $g \in
\GL_m(\RR)$ with $g_{ij}=0$ unless $(i,j) \in I$.
\end{lem}

This lemma is well known, so we only sketch the key arguments. The
commutative group $\bT^m$ acts by conjugation on the Lie algebra of $H$,
which therefore must be a direct sum of simultaneous eigenspaces of the
elements of $\bT^m$ in their conjugation action on the space of $m
\times
m$-matrices. These simultaneous eigenspaces are the one-dimensional
subspaces spanned by the $E_{ij}$, so the Lie algebra of $H$ is spanned
by some of these matrices. For this argument, see [\citet{Borel91}, Section~8.17]. The inclusion $\bT^m \subseteq H$ implies
that the $E_{ii}$ are all in the Lie algebra, and for $E_{ij},E_{jk}$
in the Lie algebra with $i \neq k$, also the commutator
$[E_{ij},E_{jk}]=E_{ik}$ lies in the Lie algebra. The earliest relation
to pre-orders that we could find is the paper \citet{malysev77}.

%
%


Next, we determine which $E_{ij}$ lie in $\ggg$.

\begin{prop}\label{propH0gens}
For $i, j\in[m], $ the matrix $E_{ij}$ lies in $\ggg$ if and only if
${j}\preccurlyeq{i}$. As a consequence, the identity component of $G$
is the group $G^0=\{g \in\GL_m(\RR) \mid g_{ij} =0 \mbox{ if } j
\not
\preccurlyeq i\}$ from the \hyperref[secintro]{Introduction}.
\end{prop}

\begin{pf}
The element $E_{ij}$ with $i \neq j$ lies in $\ggg$ if and only if
the one-parameter group $(I+tE_{ij}),\ t \in\RR$ lies in $G$, that
is, maps $\SG$ into itself. Pick $K \in\SG$ with nonzero entries on\vadjust{\goodbreak}
the diagonal and at all positions corresponding to edges of~$\cG$. We
have $(I+t E_{ij}) \cdot K = (I-t E_{ji}) K (I-tE_{ij})$---this takes
into account the inverses and the transpose in the definition of the
action. This action has the effect of subtracting $t$ times the $i$th
row of $K$ from the $j$th row and subtracting $t$ times the $i$th
column from
the $j$th column. For suitable $t$ this will create zeroes at positions
corresponding to nonedges of $\cG$ \emph{unless} the positions of the
nonzeroes in the $i$th row are among the positions of the nonzeroes in
the $j$th row. This shows that $N(i) \cup\{i\} \subseteq N(j) \cup
\{j\}$ is necessary for $E_{ij}$ to lie in $\ggg$; and repeating the
argument for general $K$ shows that it is also sufficient. The second
statement now follows from Lemma~\ref{lemifdiag}.
\end{pf}

Recall the running example $P_{3}\dvtx \stackrel{2}{\bullet}
-\stackrel{1}{\bullet}-\stackrel{3}{\bullet}$ from the
\hyperref[secintro]{Introduction}. By
Proposition~\ref{propH0gens}, the Lie algebra $\ggg$ is spanned by
$E_{11},E_{22},E_{33}$ together with $E_{21}$ and $E_{31}$. The element
$E_{21}$ lies in $\ggg^{0}$ because $N(2)\cup\{2\} =\{1,2\}\subseteq
N(1)\cup\{1\}=\{1,2,3\}$. The inverse containment does not hold, so
$E_{12}$ does not lie in $\ggg^{0}$. The group $G^{0}$ consists of
invertible matrices of the form
\[
\left[ %
\matrix{ * & 0 & 0
\vspace*{2pt}\cr
* & * & 0
\vspace*{2pt}\cr
* & 0 & * }
\right],
\]
where the asterisk denotes an element which can be nonzero.

\begin{figure}

\includegraphics{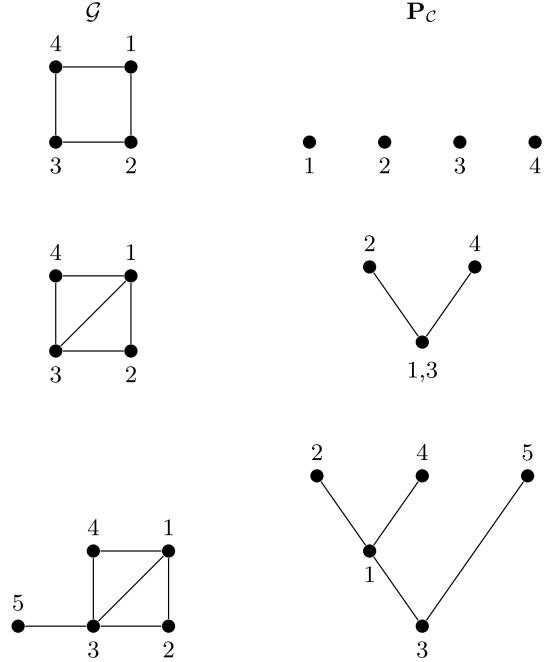}

\caption{Three graphs and Hasse diagrams of the corresponding posets
$\bP_{\cC}$.}\label{figPG}
\end{figure}

It is useful in the remainder of the paper to have a thorough
understanding of the pre-order $\preccurlyeq$. It can also be
described in terms of the collection $\cC$ of \emph{maximal cliques}
in the graph $\cG$, as follows: $j \preccurlyeq i$ if and only if
every $C
\in\cC$ containing $j$ also contains $i$. Recall that $\preccurlyeq$
determines an equivalence relation $\sim$ on $[m]$. It also determines
a partial order on $[m]/\sim$, still denoted $\preccurlyeq$, defined
by $\bar{i} \preccurlyeq\bar{j}$ if $i \preccurlyeq j$. We denote the
poset $([m]/\sim,\preccurlyeq)$ by $\bP_\cC$; it was first introduced
in~\citet{letacmassam2007} but appeared also in other related contexts
in \citet{andersson2010,drton2008}. In Figure~\ref{figPG}, we show three
graphs and the \textit{Hasse diagrams} of the
corresponding posets $\bP_{\cC}$. We note in passing that not all
posets arise as $\bP_{\cC}$ for some $\cG$. Two counterexamples
are given in Figure~\ref{fignonposets}. A more detailed study of the
structure of $\mathbf{P}_\cC$ is provided in the supplementary materials [\citet{supp}].

\begin{figure}

\includegraphics{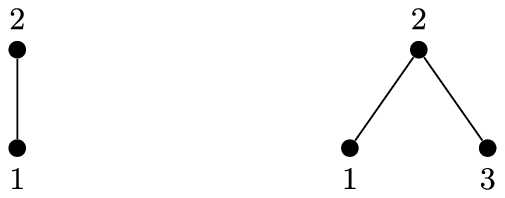}

\caption{Two posets that do not arise as $\bP_{\cC}$ for any $\cG$ with
collection $\cC$ of maximal cliques.}\label{fignonposets}
\end{figure}

\begin{rem}
Imagine relabeling the vertices of $\cG$ by $[m]$ in such a way that
the equivalence classes of $\sim$ are consecutive intervals and such
that an inequality $\bar{j}\prec\bar{i}$ between equivalence classes
implies that the interval corresponding to $\bar{j}$ contains smaller
integers than the interval corresponding to $\bar{i}$. Then the matrices
in $G^0$ are block lower triangular with square blocks along the diagonal
corresponding to the equivalence classes. From this it is easy to see
that $G^0$ is connected in the Zariski topology, but not in the ordinary
Euclidean topology. Its number of Euclidean components is $2^{|[m]/\sim|}$,
corresponding to sign patterns of the determinants of the diagonal
blocks.\vadjust{\goodbreak}
\end{rem}

\begin{rem}
The analogue of $G^{0}$ has been studied for other Gaussian models. For lattice
conditional independence models this group was named the group
of generalized block-triangular matrices with lattice structure
[see \citet{andersson1993}, Section~2.4]. The link between lattice conditional
independence models and certain Gaussian graphical models is
discussed in \citet{andersson1995} and in the supplementary materials [\citet{supp}].
\end{rem}

\subsection{The component group}

Now that we have determined the identity component $G^0$ of $G$, we set
out to describe the quotient $G/G^0$, known as the component group.
In the \hyperref[secintro]{Introduction} we observed that for our running example $P_3$
the permutation matrix
\[
\lleft[\matrix{ 1 & 0 & 0\vspace*{2pt}
\cr
0 & 0 & 1\vspace*{2pt}
\cr
0 & 1 & 0}
\rright],
\]
lies in $G$ but not in $G^0$. The key to generalizing this observation
is the following.

\begin{prop}\label{lemasperm}
Every element $g\in G$ can be written as $g=\sigma g_0$, where $g_0\in
G^0$ and $\sigma$ is a permutation matrix contained in $G$.
\end{prop}

\begin{pf}
The subgroup $H:=g^{-1} \bT^m g$ is a maximal (real, split) torus in
the real algebraic group $G^0$. By a standard result in the theory of
algebraic groups [see, e.g., \citet{Borel91}, Theorem 15.14], maximal
tori are conjugate under $G^0$. Hence, there exists a $g_0 \in G$ such
that $\bT^m=g_0^{-1} H g_0$. Then $\bT^m=(gg_0)^{-1} \bT^m (gg_0)$,
that is, $gg_0$ normalizes $\bT^m$. But the normalizer of $\bT^m$
in $\GL_m(\RR)$ consists of monomial matrices, that is, $gg_0$ equals
$\sigma t$ with $t \in\bT^m$ and $\sigma$ some permutation matrix.
Hence, $g=\sigma(tg_0^{-1})$. Here the second factor is an element of
$G^0$, so that $\sigma$ is a permutation matrix contained in $G$.
\end{pf}

We can now prove Theorem~\ref{thmmG}.

\begin{pf*}{Proof of Theorem~\ref{thmmG}}
By Proposition~\ref{lemasperm} every element of $G$ can be written as
$\sigma
g_0$ with $g_0$ an element of $G^0$ and $\sigma$ a permutation matrix
belonging to $G$, that is, preserving the zero pattern of matrices in
$\SG$.
The only such permutation matrices are those coming from automorphisms
of $\cG$. This proves that $G=\Aut(\cG) G^0$, where we identify the
automorphism group $\Aut(\cG)$ with the group of corresponding permutation
matrices. This proves the theorem.
\end{pf*}

As explained in the \hyperref[secintro]{Introduction}, the expression $G=\Aut(\cG) G^0$
is not minimal in the sense that $\Aut(\cG)$ and $G^0$ may intersect.
To get rid of that intersection, we define $\tcG$ to be the graph with
vertex set $[m]/\sim$ and an edge between $\bar{i}$ and $\bar{j}$ if there
is an edge between $i$ and $j$ in $\cG$. Define $c\dvtx [m]/\sim\to\NN,\
\bar{i} \mapsto|\bar{i}|$ and view $c$ as a coloring of the vertices
of $\tcG$ by natural numbers. Let $\Aut(\tcG,c)$ denote the group of
automorphisms of $\tcG$ preserving the coloring. There is a lifting
$\ell\dvtx  \Aut(\tcG,c) \to\Aut(G)$ defined as follows: the element
$\tau\in\Aut(\tcG,c)$ is mapped to the unique bijection $\ell(\tau)\dvtx [m]
\to[m]$ that maps each equivalence class $\bar{i}$ to the equivalence
class $\tau(\bar{i})$ by sending the $k$th smallest element of $\bar{i}$
(in the natural linear order on $[m]$) to the $k$th smallest element
of $\tau(\bar{i})$, for $k=1,\ldots,|\bar{i}|$.

\begin{figure}

\includegraphics{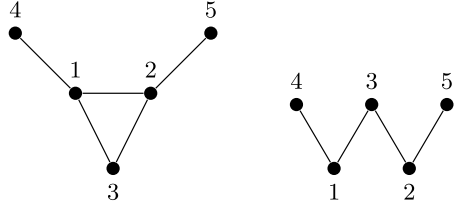}

\caption{The bull graph and the corresponding $\bP_{\cC}$ on the
right.}\label{figbadcase}
\end{figure}

\begin{thmm} \label{thmmGsemidir}
The group $G$ equals $\ell(\Aut(\tcG,c)) G^0$, and the intersection
$\ell(\Aut(\tcG,c)) \cap G^0$ is trivial, so $G$ is the semidirect product
$\ell(\Aut(\tcG,c)) \ltimes G^0$.\vadjust{\goodbreak}
\end{thmm}

\begin{pf}
By the proof of Theorem~\ref{thmmG}, any $g \in G$ can be written as
$\sigma g_0$ with $\sigma\in\Aut(\cG)$ and $g_0 \in G^0$. Since
$\sim$ is defined entirely in terms of the graph $\cG$, the graph
automorphism $\sigma$ satisfies $i \sim j \Leftrightarrow\sigma(i)
\sim\sigma(j)$. This implies that $\sigma$ determines an automorphism
$\tau\in\Aut(\tcG,c)$ mapping $\bar{i}$ to $\overline{\sigma
(i)}$. Now
$\sigma$ equals $\ell(\tau) \sigma'$ where $\sigma' \in\Aut(\cG
)$ maps
each equivalence class $\bar{i}$ into itself. But then $\sigma'$ lies in
$G^0$ and hence $g$ equals $\ell(\tau)$ times an element $\sigma' g_0$
of $G^0$. This proves the first statement. As for the second statement,
observe that a permutation matrix can have the zero pattern prescribed
by $G^0$ \emph{only} if the permutation maps each equivalence class into
itself. The only element of $\ell(\Aut(\tcG,c))$ with this property is
the identity matrix.
\end{pf}


\begin{exmp}As an example, we consider a special small graph---the
\textit{bull graph}---which is a graph on five vertices depicted in
Figure~\ref{figbadcase}. The continuous part of $G$ is given by the
poset $\mathbf{P}_\cC$ depicted on the right. There is only one
nontrivial automorphism of $\cG$. It permutes $4$ with $5$ and $1$
with~$2$. Hence, the group $G\subseteq\mathrm{GL}_5(\RR)$ consists of
matrices of the following two types:
\[
\left[ %
\matrix{ * & 0 & 0 & 0 & 0
\vspace*{2pt}\cr
0 & * & 0 & 0 & 0
\vspace*{2pt}\cr
* & * & * & 0 & 0
\vspace*{2pt}\cr
* & 0 & 0 & * & 0
\vspace*{2pt}\cr
0 & * & 0 & 0 & * }
\right] \quad\mbox{and}\quad \left[ %
\matrix{ 0 & 1 & 0 & 0 & 0
\vspace*{2pt}\cr
1 & 0 & 0 & 0 & 0
\vspace*{2pt}\cr
0 & 0 & 1 & 0 & 0
\vspace*{2pt}\cr
0 & 0 & 0 & 0 & 1
\vspace*{2pt}\cr
0 & 0 & 0 & 1 & 0 }
\right]\cdot\left[ %
\matrix{ * & 0 & 0 & 0 & 0
\vspace*{2pt}\cr
0 & * & 0 & 0 & 0
\vspace*{2pt}\cr
* & * & * & 0 & 0
\vspace*{2pt}\cr
* & 0 & 0 & * & 0
\vspace*{2pt}\cr
0 & * & 0 & 0 & *}
\right]
\]
\end{exmp}

To see Theorems~\ref{thmmG} and~\ref{thmmGsemidir} in some
further examples, see Section~\ref{secexamples}.

\begin{rem}To the coloured graph $(\widetilde\cG,c)$ we can associate
a Gaussian graphical model $M(\cG,c)$ with \textit{multivariate
nodes}, where node $\bar i$ is associated to a Gaussian vector of
dimension ${c_{\bar i}}$. This model coincides with $M(\cG)$. This
also shows, conversely, that our framework extends to general Gaussian
graphical models with multivariate nodes.
\end{rem}

\section{Existence and robustness of equivariant estimators}\label
{secequivariance}

Suppose that in the inference of the unknown concentration matrix
$K\in\SplusG$ the observed $n$-sample $\mathbf{x}\in\RR^{m \times
n}$ leads to the estimate $T(\mathbf{x})$.\vadjust{\goodbreak} Then it is reasonable
to require that the sample $g\mathbf{x}$ leads to the estimate $g
T(\mathbf{x})$. Such a map $T\dvtx \RR^{m \times n} \to\SplusG$, possibly
defined only outside some (typically $G$-stable) measure-zero set and
satisfying $T(g\bX)=gT(\bX)$ for all $g \in G$ there, is
called a ($G$-)\emph{equivariant estimator}. In this section we determine
a sharp lower bound on $n$ for an equivariant estimator to exist, and
then, building on theory from \citet{daviesgather2005}, we determine a
bound on the robustness of such estimators.

\begin{figure}

\includegraphics{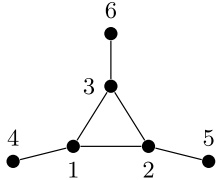}

\caption{For this graph a $G$-equivariant map exists as soon as $n\geq2$.
However, the MLE exists only when $n\geq3$.}\label{figbibull}
\end{figure}

The MLE, when it exists, is automatically $G$-equivariant, since the
likelihood function is $G$-invariant. A \emph{necessary} condition for the
MLE to exist with probability $1$ is that the sample size $n$ be at least
the largest clique size $q=\max_{C \in\cC} |C|$. A~\emph{sufficient}
condition is that $n$ be at least the maximal clique size $q^*$ in
a \emph{decomposable cover} of $\cG$, that is, a graph $\cG^*=([m],E^*)$
with $E^* \supseteq E$ that does not have induced $k$-cycles for $k
\geq
4$. The exact minimal value of $n$ for which MLE exists is not known
explicitly in general, but interesting classes of graphs were analyzed
in \citet{buhl93,barret96} and \citet{uhler2010}.

Our Theorem~\ref{thGeqexist} states that equivariant estimators of
the concentration matrix (or, equivalently by taking inverses, of
the covariance matrix) exist if and only if $n \geq\max_{i \in[m]}
|\down i|$. Note that this is weaker than the necessary condition $n
\geq q$
for the existence of MLE. Indeed, any down set $\down i$ is in fact
a clique, because $j,k \in\downb i$ implies that $j \in N(i) \cup
\{i\} \subseteq N(k) \cup\{k\}$, that is, $j$ and $k$ are either equal or
connected by an edge. The inequality $\max_{i \in[m]} |\down i| \leq
\max_{C \in\mathcal{C}} |C|$ can be strict. For example, in the
graph of
Figure~\ref{figbibull} the biggest maximal clique has cardinality $3$, while
$\max_{i\in[m]}|\down i|=2$. In consequence, our result
does not shed new light on the existence of MLE; however, it does provide
necessary and sufficient conditions in the search for other equivariant
estimators.

\subsection{Existence of equivariant estimators}

We now prepare the proof of Theorem \ref{thGeqexist}. In our arguing,
we borrow some terminology from algebraic geometry: we say that some
property holds for \textit{generic} $n$-tuples $\mathbf{x}\in\RR
^{m\times
n}$ if it holds for $\mathbf{x}$ outside the zero set of some nonzero
polynomial. Note that if a property holds for generic $n$-tuples,
then it holds with probability one for the random sample $\bX$
drawn from any nondegenerate probability distribution with continuous
density function on
$\RR^{m\times n}$.

\begin{thmm}\label{thdelta0}
The minimal number $n$ for which the stabilizer in $G^{0}$ of a generic
$n$-sample $\mathbf{x}\in\RR^{m \times n}$ consists entirely of
determinant-$(\pm1)$ matrices equals $n=\max_{i\in[m]} |\down i|$. For
that value of $n$ the stabilizer of a generic $n$-sample is, in fact,
the trivial group $\{I\}$.
\end{thmm}

\begin{pf}
The condition that $g \in G^{0}$ fixes one vector
$x=(x_1,\ldots,x_m)^T \in\RR^m$ translates into $m$ linear conditions
on the entries of $g$, namely:
\[
\sum_{j\preccurlyeq i} g_{ij} x_j =
x_i\qquad \mbox{for } i=1,\ldots,m.
\]
The $i$th condition concerns only the entries in the $i$th row of
$g$. We therefore concentrate on that single row of $g$, and regard
the entries $g_{ij},\ j \preccurlyeq i$ as variables to be solved from
the linear equations above as $x$ ranges through the given $n$-tuple
$\mathbf{x}$
of vectors. Since the given $n$-tuple is generic, those equations
are linearly independent as long as $n$ is at most the cardinality
of $\down i$. Hence, they determine the $i$th row uniquely as
soon as $n$ is at least that number. Hence, as soon as $n$ is at least
the maximal cardinality of the sets $\down i$ over all $i$ the
stabilizer of a generic sample $\mathbf{x}$ will be trivial.

What remains to be checked, is that for smaller $n$ the stabilizer
of a generic $\mathbf{x}\in\RR^{m \times n}$ does \emph{not} consist
entirely of determinant-$(\pm1)$ matrices. This is most easily seen by
considering the Lie algebra of that stabilizer, which is the
set of matrices $A$ in the Lie algebra of $G^0$ satisfying the linear
conditions $A\mathbf{x}=0$. Let $i$ be a row index for which
$\down i$ has more than $n$ elements. Then the linear conditions on $A$
do not fix the $i$th row of $A$ uniquely. Moreover, by genericity, they
do not fix the diagonal entry $A_{ii}$ uniquely, either. As a consequence,
they do not determine the trace of $A$ uniquely. This shows that the
Lie algebra of the stabilizer is not contained in the Lie
algebra of trace-zero matrices. But then the stabilizer is
not contained in the Lie group of determinant-one matrices (whose Lie
algebra consists of the trace-zero matrices).
\end{pf}

\begin{pf*}{Proof of Theorem \ref{thGeqexist}, necessity of $n \geq
\max_i|\down i|$}
Assume that there exists a $G$-equivariant estimator $T\dvtx \RR^{m
\times n} \to\SplusG$, possibly defined outside some measure-zero
set. In particular, the $G^{0}$-equivariance of $T$ implies that the
$G^{0}$-stabilizer of a generic sample $\mathbf{x}$ is contained in the
$G^{0}$-stabilizer of $T(\mathbf{x})$:
\[
G^0_{\mathbf{x}} \leq G^0_{T(\mathbf{x})}.
\]
Now since $T(\mathbf{x})\in\SplusG$, the stabilizer on the
right-hand side
is a generalized orthogonal group, and hence in particular compact in the
Eculidean topology. Hence,
the stabilizer on the left-hand side must be compact, as well.
However, by (the proof of) Theorem \ref{thdelta0}, that stabilizer
is the intersection of $\GL_m(\RR)$ with an affine
subspace of $\RR^{m \times m}$. Such a set is not compact in the Euclidean
topology unless it consists of a single matrix, and this happens only
when $n\geq\max_{i\in[m]}|\down i|$.
\end{pf*}

To prove that $n \geq\max_i |\down i|$ is also \emph{sufficient}
for the existence of a $G$-equivariant estimator, we introduce the
following construction. Fix a natural number $n \geq\max_i |\down i|$
and construct a function $f\dvtx [m] \to[n]$ by induction, as follows: if $f$
has been defined on all elements of $\down i \setminus\bar{i}$, then define
$f$ on elements of $\bar{i}$ to be the increasing bijection from $\bar
{i}$ (with
the natural linear order coming from $[m]$) to the $|\bar{i}|$ smallest
elements of the set $[n] \setminus f(\downb i \setminus\bar{i})$. This
automatically guarantees that $f$ is injective on any down set $ \down
i$ and that $f\circ g=f$ for all $g \in\ell(\Aut(\tcG,c))$. Now let
$L \subseteq
\RR^{m \times n}$ be affine space of all matrices $\mathbf{x}$ with
the property that first, the matrix $\mathbf{x}[\bar{i},f(\bar{i})]$ obtained
by taking the rows labeled by $\bar{i}$ and the columns labeled by
$f(\bar{i})$ is an identity matrix for each $\bar{i} \in[m]/\sim$, and
second, the matrices $\mathbf{x}[\bar{j},f(\bar{i})]$ are zero for
all $\bar{j} \prec\bar{i}$.

In our running example $P_3$, if the sample size $n$ is at least $2$,
then $f$ maps $1$ to $1$ and $2,3$ both to $2$. The affine space $L$
then consists of all matrices of the form
\[
\lleft[\matrix{ 1 & * & * & \cdots& *\vspace*{2pt}
\cr
0 & 1 & * & \cdots& *
\vspace*{2pt}
\cr
0 & 1 & * & \cdots& *} \rright].
\]

\begin{lem} \label{lemintoL}
For generic $\mathbf{x}\in\RR^{m \times n}$, there exists a unique
$g \in G^0$
such that $g \mathbf{x}\in L$.
\end{lem}

The geometric content of this lemma is that $L$ is a \emph{slice}
transverse to (most of) the orbits of $G^0$ on $\RR^{m \times n}$. In
our running example $P_3$, one goes from a generic sample to a sample
in $L$ by first multiplying the first row by $x_{11}^{-1}$ so as to
create a one at position $(1,1)$; then subtracting a multiple of the
(new) first row from the second to create a zero at position $(2,1)$
and multiplying the second row by a constant to create a one at position
$(2,2)$; and then similarly (and independently) for the third row. All
of these operations are realized by elements of $G^0$. The following
proof in the general case is a straightforward generalization of this.

\begin{pf*}{Proof of Lemma \ref{lemintoL}}
For the existence of such a $g$, proceed by induction. Assume that the
submatrix $\mathbf{x}[\downb i \setminus\bar{i},f(\downb i \setminus
\bar{i})]$
already has the required shape, and decompose $\mathbf{x}[\downb
i,f(\downb i)]$
into blocks as follows:
\[
\lleft[\matrix{ \mathbf{x}\bigl[\downb i \setminus\bar{i},f(\downb i
\setminus\bar{i})\bigr] & \mathbf{x}\bigl[\downb i \setminus\bar{i},f(\bar{i})
\bigr]\vspace*{2pt}
\cr
\mathbf{x}\bigl[\bar{i},f(\downb i \setminus\bar{i})\bigr]
& \mathbf{x}\bigl[\bar{i},f(\bar{i})\bigr]} \rright].
\]
Then take the block matrix $g \in G^0$ which is the identity outside the
blocks labeled by $\down i \times\down i$, and which in those blocks
looks like
\[
\lleft[\matrix{ I & 0
\vspace*{2pt}\cr
g[\bar{i},\downb i \setminus\bar{i}] & g[\bar{i},\bar{i}] } \rright].
\]
Now straightforward linear algebra shows that, under the condition that
both \mbox{$\mathbf{x}[\downb i,f(\downb i)]$} and $\mathbf{x}[\downb i
\setminus
\bar{i},f(\downb i
\setminus\bar{i})]$ are full rank, there are unique choices for the as
yet unspecified components of $g$ such that
$(g\mathbf{x})[\bar{i},f(\downb i
\setminus\bar{i})]=0$ and $(g\mathbf{x})[\bar{i},f(\bar{i})]=I$.
This shows
the existence of $g$ such that $g\mathbf{x}\in L$. Uniqueness can be
proved by a similar induction.
\end{pf*}

\begin{pf*}{Proof of Theorem \ref{thGeqexist}, sufficiency of $n \geq
\max_i|\down i|$}
Now we show that $n\geq\max_{i\in[m]}|\down i|$ is also a sufficient
condition for the existence of an equivariant map $T\dvtx \RR^{m
\times n} \to\SplusG$, defined for generic
samples $\mathbf{x}$. Indeed, by construction, the space $L$ is stable under
$\ell(\Aut(\tcG,c))$. Fix any $\ell(\Aut(\tcG,c))$-equivariant
map $T\dvtx L \to
\SplusG$. Such maps exist and can be found as follows: take $T'\dvtx L \to
\SplusG$ \emph{any} map, and then define
\[
T(\bX):=\frac{1}{|\Aut(\tcG,c)|} \sum_{g \in\ell(\Aut(\tcG,c))} g \cdot
T'\bigl(g^{-1} \bX\bigr),
\]
an average over the finite group $\Aut(\tcG,c)$.

We claim that $T$ extends to a unique $G$-equivariant map $\RR^{m
\times n} \to\SplusG$ defined almost everywhere. Indeed, this
extension is defined as follows: given a generic sample $\mathbf{x}$,
find the unique $g_0 \in G^0$ such that $g_0\mathbf{x}\in L$, and set
$T(\mathbf{x}):=g_0^{-1} \cdot T(g_0 \mathbf{x})$. Checking that the map
$T$ thus defined (almost) everywhere is both $\Aut(\tcG,c)$-equivariant
and $G^0$-equivariant is straightforward. This proves the existence
part of Theorem~\ref{thGeqexist}.
\end{pf*}

\begin{rem}
We stress that, apart from giving necessary and sufficient conditions
for the existence of a $G$-equivariant estimator, the proof of Theorem
\ref{thGeqexist} actually yields the general structure of \emph{any}
such estimator. Of course, the usefulness (bias, robustness, etc.) of
an equivariant estimator thus constructed depends on the (free) choice
of $T'$, that is, on the restriction of $T$ to the slice $L$. We do not
know at present good conditions on $T'$ that ensure usefulness of $T$.
\end{rem}

\begin{rem}
Note that the maps $T\dvtx \RR^{m \times n} \to\SplusG$ constructed in
the proof of Theorem~\ref{thGeqexist} are merely $G$-equivariant, and
not necessarily invariant under permutation of the sample points. It
is easy to see, though, that the lower bound $n \geq\max_i|\down i|$
also implies the existence of $G$-equivariant estimators that \emph{are}
invariant under permutations of the sample points. Indeed, simply replace
$T$ by its group average $\bX\mapsto\frac{1}{n!} \sum_{\sigma
\in S_n} T(\bX^\sigma)$.
\end{rem}

\subsection{Robustness}

An important notion for the robustness of parameter estimators is that
of breakdown points [\citet{hampel71,donohohuber}]. In a simple
univariate situation, if the estimator
is given by the sample mean, then a (large) change made to one of the
observations leads to an arbitrarily large change in the value of the
estimator. On the other hand, if the estimator is the sample median, then
changing one observation in a sample of size larger than two cannot lead
to arbitrarily large changes in the estimator. This feature makes the
median more robust to outliers in the sample. The (\textit{finite sample})
\textit{breakdown point} of an estimator $T$ at an $n$-sample $\bX=\mathbf
{x}$ is
the minimal number of components of $\mathbf{x}$ that need to be altered
to force arbitrarily large changes in the value of the estimator; this
quantity is usually normalized by the sample size $n$. For example,
the sample mean above has breakdown point $1/n$ while the sample median
has breakdown point roughly $1/2$ (in fact, both independently
of $\mathbf{x}$). So when it comes to robustness, the estimator with the
highest breakdown point is preferred.

In the multivariate Gaussian setting, when estimating the
concentration matrix (or the covariance matrix), the change in the
estimator value is often measured by means of the \emph{pseudo-metric}
$D$ on $\cS^+_m$ [see, e.g., \citet{daviesgather2005}]
\[
D(K_{1},K_{2}) = \bigl|\log\det\bigl(K_{1}K_{2}^{-1}
\bigr)\bigr|.
\]
For graphical models, robustness issues have been rarely looked at so far,
although it has been known for some time that the classical estimators
and model selection procedures are vulnerable to contaminated
data [\citet{gottard2007,kuhnt2003}]. First, approaches toward robust
covariance estimators for undirected Gaussian graphical models can
be found in \citet{becker2005iterative,gottard2010robust}. These
papers suggest to replace the sample covariance matrix
by the reweighted minimum covariance determinant (MCD) estimator. The
paper \citet{miyamura2006robust}
proposes an M-type estimator instead. Both in \citet{finegold2011robust}
and in \citet{Vogel01122011} the assumption of normality is
discarded, and replaced by the $t$-distribution or the general
elliptical distribution, respectively, to model heavy tails.

Our modest contribution to robustness issues is an upper bound on the
finite sample breakdown point for $G$-equivariant estimators of the
concentration matrix for the graphical model $M(\cG)$. To
this end,
we specialize one of the key ideas from~\citeauthor{daviesgather2005}
(\citeyear
{daviesgather2005,davies2007breakdown}) to
our setting. Suppose we have an $n$-sample $\mathbf{x}\in\RR^{m
\times n}$
and an equivariant estimator $T\dvtx \RR^{m \times n} \to\SplusG$ of the
concentration matrix. Assume that there exists an element $g \in G$
with $|\det
g| \neq1$ that fixes (at least) $k$ of the $n$ sample
points $\mathbf{x}_1,\ldots, \mathbf{x}_k$ of the sample $\mathbf
{x}$. Define $d=\lceil\frac{n-k}{2}\rceil$ and let
\[
\by= \bigl(\mathbf{x}_1,\ldots,\mathbf{x}_k,\ldots,
\mathbf {x}_{n-d},g^l \mathbf{x}_{n-d+1},
\ldots,g^l \mathbf{x}_n\bigr).
\]
Since $k+d\geq n-d$, for each natural number $l$ both $\by$ and
$g^{-l} \by$ contain at least $n-d$ points of the original sample
$\mathbf{x}$. By the triangle inequality, we have $D(T(\by
),T(g^{-l}\by))\leq D(T(\mathbf{x}),T(\by))+D(T(\mathbf
{x}),T(g^{-l}\by))$
and on the other hand
\begin{eqnarray*}
D\bigl(T( \by),T\bigl(g^{-l}\by\bigr)\bigr)&=&D\bigl( T(\by),
\bigl(g^{T}\bigr)^l T(\by)g^l\bigr)
\\
&=&\bigl|\log\det\bigl(\bigl(g^{T}\bigr)^l T(\by)
g^{l} T(\by)^{-1}\bigr)\bigr| =l\bigl|\log\bigl(\det g^2
\bigr)\bigr|,
\end{eqnarray*}
which is unbounded as $l \to\infty$. Hence, changing not more than
$d=\lceil\frac{n-k}{2}\rceil$
of the sample points in $\mathbf{x}$ can already
lead to arbitrarily large changes in the estimated concentration matrix,
so that the finite sample breakdown point of $T$ at $\mathbf{x}$ is at
most $d/n$. We now state and prove our upper bound on the
robustness of equivariant estimators at generic samples.

\begin{prop} \label{proprobustness}
Assume that $n \geq\max_i |\down i|$. Then for any
$G$-equivariant
estimator $T\dvtx \RR^{m \times n} \to\SplusG$ the finite sample breakdown
point at a generic sample $\mathbf{x}$ is at most $\lceil(n-\max_i
|\down i|+1)/2\rceil/n$.
\end{prop}

\begin{pf}
By Theorem~\ref{thdelta0}, there exist matrices $g \in G^0$ with
determinant $\neq\pm1$ that fix the first $k=\max_i|\down i| - 1$ sample
points. Now the proposition follows from the discussion preceding it.
\end{pf}

\begin{rem}
Writing $q:=\max_i |\down i|$, note that $q \leq m$ with equality if
and only if $\cG$ is the complete graph, and that $q \geq1$, with
equality if and only if for each edge $\{i,j\} \in E$ the vertex $i$
has neighbors that are not connected to $j$ (and vice versa). Examples
of such graphs are $m$-cycles with $m \geq4$. Trees with $m \geq3$
vertices are examples of graphs with $q=2$.

Note that for graphs with small $q$ the upper bound in
Proposition~\ref{proprobustness} is close to $1/2$, even for relatively
small sample sizes $n$. On the other hand, the MLE, as pointed out for
example in
\citet{robuststatistics}, is typically
the least robust with respect to potential outliers in the sample space.
Although we do not know whether the upper bound in the proposition is
attained for any sensible estimator, our results do suggests the quest
for more robust estimators, especially for graphs with small $q$.
\end{rem}


\section{The maximal invariant}\label{secorbequiv}

In this section, we discuss a $G^0$-invariant map $\tau$ on the
space $\RR^{m \times n}$ of $n$-samples, defined almost everywhere,
and prove that it is maximal in the sense that for two samples
$\mathbf{x},\by$ in the domain of definition of $\tau$ the equality
$\tau(\mathbf{x})=\tau(\by)$ implies that $\mathbf{x},\by$
are in the same $G^0$-orbit.

Recall from the \hyperref[secintro]{Introduction} that $\tau$ is defined as
\[
\tau\dvtx \mathbf{x}\mapsto \bigl(\mathbf{x}[\downb i]^T \bigl(
\mathbf{x}[\downb i] \mathbf {x}[\downb i]^T\bigr)^{-1}
\mathbf{x}[\downb i] \bigr)_{\bar{i}\in[m]/\sim},
\]
where we assume from now on that $n$ is at least $|\down
i|$, and where
$\tau$ is defined on $n$-samples where $\mathbf{x}[\downb i]$ has
full rank for
all $i$. Before we proceed, we recall the following known lemma.

\begin{lem}
Let $k \leq n$ be natural numbers, and consider the action of $\GL
_k(\RR)$
on $\RR^{k \times n}$. Let $U$ be the open subset of the latter space
consisting of matrices of full rank $k$. Then the map $\varphi\dvtx  U
\mapsto
\RR^{n \times n}$ mapping $\mathbf{x}$ to $\mathbf{x}^T (\mathbf
{x}\mathbf{x}^T)^{-1} \mathbf{x}$ is a
maximal invariant for the action of $\GL_k(\RR)$ on $U$.
\end{lem}

\begin{pf}
First, to see that $\varphi$ is $\GL_k(\RR)$-invariant, compute
\[
\varphi(g\mathbf{x})=\mathbf{x}^T g^T \bigl(g \mathbf{x}
\mathbf{x}^T g^T\bigr)^{-1} g\mathbf{x}=
\varphi(x).
\]
Second, to see that $\varphi$ is maximal, note that the row space of
$\mathbf{x}\in U$
is also the row space of $\varphi(\mathbf{x})$. Hence, if
$\varphi(\by)=\varphi(\mathbf{x})$ for
a second $\by\in U$, then $\by$ has the same row space as $\mathbf
{x}$. But this
means that there exists a $g \in\GL_k(\RR)$ with $g \mathbf{x}=\by$.
\end{pf}

The proof of the lemma shows that $\varphi(\mathbf{x})$ determines
the row
space of $\mathbf{x}$ (and is determined by that!). Now we can prove
Theorem~\ref{thmmmaxinv}, which states that $\tau$ is a maximal
$G^0$-invariant. This generalizes Example 6.2.3 in \citet{lehmannromano},
which deals with the case of complete graphs.

\begin{pf*}{Proof of Theorem~\ref{thmmmaxinv}}
The $G^0$-invariance of each of the components of $\tau$ follows from
the observation that $(g \mathbf{x})[\downb i]=g[\downb i] \mathbf
{x}[\downb i]$, together
with the computation in the proof of the preceding lemma.

For maximality, assume that $\tau(\mathbf{x})=\tau(\by)$. This
means that
the row space of $\mathbf{x}[\downb i]$ equals that of $\by[\downb
i]$, for
all $i$. If, by induction, we have replaced $\mathbf{x}$ by an element
in its
orbit and achieved that $\mathbf{x}[\downb i \setminus\bar{i}]=\by
[\downb i
\setminus\bar{i}]$, then it follows that $\by[\bar{i}]=h_1 \mathbf
{x}[\bar{i}]
+ h_2 \mathbf{x}[\downb i \setminus\bar{i}]$ for a suitable invertible
$\bar{i}\times\bar{i}$-matrix $h_1$ and a suitable full-rank $\bar{i}
\times(\down i \setminus\bar{i})$-matrix $h_2$. These matrices $h_1,h_2$
can be assembled into a block matrix $g_0 \in G^0$ (as in the proof
of Lemma~\ref{lemintoL}) such that $(g_0 \mathbf{x})$ coincides with
$\mathbf{x}$
outside the $\bar{i}$-labeled rows and with $\by$ in the $\bar{i}$-labeled
rows. Doing this for all equivalence classes $\bar{i}$ from the bottom
to the top of $\bP_\cC$, we move $\mathbf{x}$ to $\by$ by an
element of $G^0$.
\end{pf*}

Since every invariant test depends on $\mathbf{x}$ only through the
value of the maximal invariant [\citet{lehmannromano}, Section~6.2],
Theorem~\ref{thmmmaxinv} paves the way for $G^0$-invariant tests, for
example, for
testing the hypothesis that the distribution of $X$ lies in $M(\cG)$
against the null-hypothesis that it does not. A more general question
involves testing two alternative (typically nested) graphical models
corresponding to graphs $\cG_1,\cG_2$ on $[m]$. For this, it is natural
to develop tests that are invariant with respect to matrices stabilizing
\emph{both} models. The identity component of the group of such matrices
consists of all $g$ with $g_{ij}=0$ unless $j \preccurlyeq i$ in \emph{both} pre-orders coming from $\cG_1,\cG_2$. The same construction as
above, now applied to the intersection of the pre-orders, gives the
maximal invariant for this group. A simple example of a $G^0$-invariant
test is the deviance test [see, e.g., \citet{lauritzen96}, Section~5.2.2].


\section{Equivariance in the transitive case} \label{sectransitive}

When $G$ acts transitively on
$\cS^+_\cG$ then $M(\cG)$ forms an exponential transformation family
[see \citet{nielsentransformation82}], which gives very efficient
tools for dealing with the ancillary statistics in the hypothesis
testing and inference. In particular the $p^*$-formula of \citet
{barndorffMLE}, which gives an approximation for the density of the
maximum likelihood estimator given the ancillary statistics is exact
and the ancillary statistics is given by the maximal invariant $\tau
(\bX)$.

The following result tells us when the graphical Gaussian
model $M(\cG)$ is an exponential transformation family under
the group $G$ (cf. Theorem \ref{thmmorbitspacedim}).

\begin{thmm}[{[Theorem 2.2, \citet{letacmassam2007}]}]\label{thmmmassam} Let
$\cG=([m],E)$ be an undirected graph. Then $G$ acts transitively
on $M(\cG)$ if and only if one of the following equivalent conditions holds:
\begin{itemize}
\item for any two neighbors $i,j\in
[m]$ we have either $i\preccurlyeq j$ or $j\preccurlyeq i$;
\item$\cG$ is decomposable and does not contain a $4$-chain
${\bullet}-{\bullet}-{\bullet}-{\bullet}$ as an induced subgraph;
\item the Hasse diagram of $\mathbf{P}_\cC$ is a tree with a unique minimum.
\end{itemize}
\end{thmm}

As we show in the supplementary materials [\citet{supp}], the transitive case is
precisely the case when $M(\cG)$ corresponds to a lattice conditional
independence model. We also prove there the following lemma.

\begin{lem}\label{lemequconds}
If $\cG$ satisfies the conditions of Theorem \ref{thmmmassam}, then
$\max_i |\down i|$ is equal to the size of the biggest maximal clique
of $\cG$. In particular a $G$-equivariant estimator exists with
probability one if and only
is the MLE estimator exists with probability one.
\end{lem}


In the transitive case construction of a $G$-equivariant estimator
is particularly straightforward [\citet{eaton1989}, Chapter~6, Example
6.2]. This generalizes the
case of a star-shape graph analyzed in \citet{sun2005estimation}. Let
$\widehat{\Sigma}$ be the MLE of the covariance matrix and define
$S(\bX)=n\widehat{\Sigma}$. Because $S(\bX)$ is a sufficient statistic,
without loss we can assume that every estimator based on the full
sample satisfies
$T(\bX)=T(S(\bX))$. Since $S(\bX)^{-1}\in
\cS^+_\cG$, there exists $h\dvtx \RR^{m\times n}\to G^0$ such that
$S^{-1}(\bX)=h(\bX)^T h(\bX)$. The construction of $h$ follows by
the fact that in the transitive case there exists a well defined map
$\phi\dvtx  \cS^+_\cG\to G/G_I$, where $G_I$ is the stabiliser of the
identity matrix. This map is the inverse of the canonical map from
$G/G_I$ to $\cS^+_\cG$. Then $h$ is just a composition of $S^{-1}(\bX
)\dvtx \RR^{m\times n}\to\cS^+_\cG$ followed by $\phi$. By
$G$-equivariance,
\[
T(\bX)=T\bigl(S(\bX)\bigr)=T\bigl(h(\bX)^Th(\bX)\bigr)=h(
\bX)^T T(I)h(\bX),
\]
where $T(I)\in\cS^+_\cG$. We have just shown the following result.
%
\begin{prop}Let $\cG$ be a decomposable graph without induced
$4$-chains. Define $S(\bX)=n\widehat{\Sigma}$ as above. Then every
$G$-equivariant estimator of the concentration matrix is of the form
\[
T(\bX) = \bigl(h_0 h(\bX)\bigr)^T h_0h(\bX),
\]
where $h_0\in G^0$ is a constant matrix and $h\dvtx \RR^{m\times n}\to
G^0$ is such
that $S(\bX)=h(\bX)^Th(\bX)$.
\end{prop}

Since the function $h$ is uniquely identified the only way to obtain
different equivariant estimators is by varying the constant matrix
$h_0$. This can be done with different optimality criteria in mind. An
interesting problem is to find $h_0$ such that $T^{-1}$ is an unbiased
estimator of the concentration matrix. Another motivation is that the
MLE for lattice conditional independence models (and hence for $M(\cG
)$ in the transitive case by the theorem in the supplementary
materials [\citet{supp}]) is not admissible [see \citet{konno2001inadmissibility}]. A
relevant question is to analyse equivariant estimators minimizing risk
related to certain loss functions. This analysis has been already done
for star-shaped models by \citet{sun2005estimation}.

\section{Orbits of $G$ on $\SplusG$} \label{secorbitS}

Given an undirected graph $\cG=([m],E)$, we have determined the group $G
\subseteq\GL_m(\RR)$ of all invertible linear maps $\RR^m \to\RR^m$
stabilizing the cone $\SplusG$. Theorem \ref{thmmmassam} characterizes
when $M(\cG)$ is a transformation family, that is, when $G$ has a single
orbit on $\SplusG$. For general $\cG$, the orbit space $\SplusG/G$---like
many quotients of manifolds by group actions---can conceivably be very
complicated. In this section, we compute its first natural invariant,
namely, its dimension. In the zero-dimensional case, we recover the class
from Theorem~\ref{thmmmassam}.

Basic Lie group theory tells us that $\dim\SplusG/G$ equals
\[
\dim\SplusG- \dim G + \dim G_K,
\]
where $G_K$ is the stabilizer of a generic concentration matrix in
$\SplusG$. In this expression, the first term equals $m+|E|$ and the second
term equals $\dim G^0 = \sum_{\bar{i} \in[m]/\sim} |\bar{i}| \cdot
|\down i|$,
so it suffices to determine the generic stabilizer dimension. Note that
for the dimension it does not matter whether we consider the stabilizer in
$G$ or in $G^0$. The following theorem makes use of the colored quotient
graph $(\tcG,c)$ from Section~\ref{secG}.

\begin{prop}\label{propdimstab}
The dimension of the stabilizer $G^0_K$ in $G^0$ of a generic matrix in
$\SplusG$ equals $\sum_{\bar{i} \in[m]/\sim} {n_{\bar{i}}\choose2}$,
where $n_{\bar{i}}$ is defined by
\[
n_{\bar{i}}:=\max \biggl\{0, |\bar{i}| - \biggl( \sum
_{\bar{j}\in
N(\bar{i}), \bar{i} \not\preccurlyeq\bar{j} \not
\preccurlyeq\bar{i}} |\bar{j}| \biggr) \biggr\},
\]
where the sum ranges over all neighbors $\bar{j}$ of $\bar{i}$ in
$\tcG$
that are not comparable to $\bar{i}$ in the partial order
$\preccurlyeq$.
\end{prop}

In words: starting from $\tcG$, one deletes all edges between vertices
that are comparable in the partial order $\preccurlyeq$, and one subtracts
from $|\bar{i}|$ the sum of the $|\bar{j}|$ for all neighboring $\bar{j}$
in the new graph. If the result is positive, then this is $n_{\bar{i}}$;
otherwise, $n_{\bar{i}}$ is zero.

The expression above suggests that the identity component of $G^0_K$
is a product of special orthogonal groups of spaces of dimensions
$n_i$, which is indeed what the proof of this proposition, given in
the supplementary materials [\citet{supp}], will show. We now use the proposition to
explain the combinatorial procedure in the \hyperref[secintro]{Introduction}.

\begin{pf*}{Proof of Theorem~\ref{thmmorbitspacedim}}
By Proposition~\ref{propdimstab}, we need to compute
\[
\bigl(m + |E|\bigr) - \sum_{\bar{i} \in[m]/\sim} |\bar{i}||\down i| + \sum
_{\bar{i} \in[m]/\sim} \pmatrix{n_{\bar{i}}
\cr
2}.
\]
The term $m$ cancels against the diagonal entries in $G^0$ in the second
term. Recall that in Theorem~\ref{thmmorbitspacedim} we colored an edge
$\{i,j\}$ in $\cG$ blue, green or red according to whether zero, one,
or two of the statements $i \preccurlyeq j$ and $j \preccurlyeq i$ hold.
The term $|E|$ counts blue plus green plus red. What remains of the second
term after cancelling the diagonal entries against $m$ counts green
edges once
and red edges twice. Thus, the first two terms count blue edges minus
red edges. Finally, the last term counts the number of red edges that
survive when blue edges are deleted one by one.
\end{pf*}

We conclude with few examples of the use of Proposition~\ref{propdimstab}.

\begin{exmp}
Let $\tcG$ be the bull graph in Figure~\ref{figbadcase}, with each vertex
representing an equivalence class in $[m]/\sim$ with cardinality $c_i$
for $i=1,\ldots,5$. In this case the only pair of connected but not
comparable vertices is $(1,2)$. With the convention that ${m\choose2}=0$
if $m\leq0$, Proposition~\ref{propdimstab} shows that the dimension
of the
stabilizer of a generic matrix in $\SplusG$ is
\[
\pmatrix{c_1-c_2\cr 2}+\pmatrix{c_2-c_1
\cr 2}+\pmatrix{c_3\cr 2}+\pmatrix{c_4\cr 2}+\pmatrix{c_5
\cr 2}.
\]
\end{exmp}

\begin{exmp}Let $\widetilde{\cG}$ be a tree, where each vertex $v$
represents an equivalence class with cardinality $c_v$. In this case,
the dimension of the stabilizer of a generic matrix in $\cS^+_\cG$ is
\[
\sum_{(u,v)\in\mathrm{inner}} \biggl(\pmatrix{c_u-c_v
\cr 2}+\pmatrix{c_v-c_u\cr 2} \biggr)+\sum
_{i\in\mathrm{leaves}}\pmatrix{c_i\cr 2},
\]
where the first sum
is over all the inner vertices of $\widetilde{\cG}$ and the second sum
is over all the leaves (vertices of valency $1$) of $\widetilde{\cG
}$. In
particular, if for some $c$ we have that $c_i=c$ for all $i\in C$ then this
formula degenerates to $l{c\choose2}$, where $l$ is the number of leaves.
\end{exmp}

\section{Small examples}\label{secexamples}


Let $S_{m}$ denote the symmetric group on $[m]$, $D_{m}$ the dihedral
group of graph isomorphisms of an $m$-cycle. Also recall that $\mathbf
{T}^{k}\simeq
(\GL_{1}(\RR))^{k}$ denotes the group of all diagonal invertible
$k\times
k$ matrices. In Table~\ref{tabsmall}, we provide the full description
of $G$ for all undirected graphs on $m=2,3,4$ vertices.

\begin{table}
\caption{Small undirected graphs $\cG$, corresponding groups $G^0$ and
$\mathrm{Aut}(\widetilde{\cG},c)$ up to isomorphism}\label{tabsmall}
\begin{tabular*}{\textwidth}{@{\extracolsep{\fill}}lccc@{}}
\hline
\multicolumn{1}{@{}l}{$\bolds{\cG}$} & \multicolumn{1}{c}{$\bolds{\bP_{\cC}}$} &
\multicolumn{1}{c}{$\bolds{G^{0}}$} & \multicolumn{1}{c@{}}{$\bolds{\mathrm{Aut}(\widetilde{\cG},c)}$} \\
\hline

\includegraphics{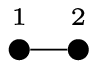}
 & \includegraphics{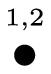} & $\GL_{2}(\RR)$ &
$\{\mathrm{id}\}$\\[3pt]

\includegraphics{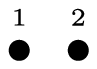}
 &\includegraphics{1130i02}& $\mathbf{T}^{2}$&$S_{2}$\\[6pt]

\includegraphics{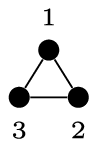}
 &
\includegraphics{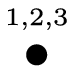}
 & $\GL_{3}(\RR)$& $\{\mathrm{id}\}$\\[3pt]

\includegraphics{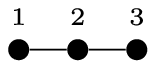}
 &
\includegraphics{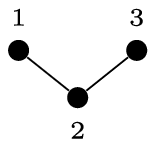}
 & $\scriptsize\left[\matrix{
 * & * & 0\vspace*{2pt}\cr
 0 & * & 0\vspace*{2pt}\cr
 0 & * & *}\right]$& $S_{2}$\\[16pt]

\includegraphics{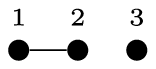}
 & \includegraphics{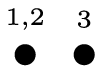} & $\GL_{2}(\RR)\times\mathbf{T}^{1}$& $\{\mathrm{id}
\}$\\[3pt]

\includegraphics{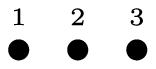}
 & \includegraphics{1130i06} & $ \mathbf{T}^{3}$& $ S_{3}$\\[6pt]

\includegraphics{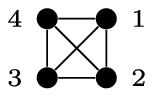}
 & \includegraphics{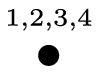} & $\GL_{4}(\RR)$& $\{\mathrm{id}\}$ \\[3pt]

\includegraphics{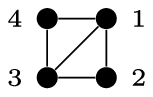}
 & \includegraphics[raise=-17pt]{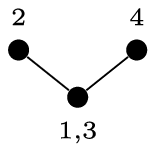} & $\scriptsize\left[\matrix{
 * & 0 & * & 0\vspace*{2pt}\cr
 * & * & * & 0\vspace*{2pt}\cr
 * & 0 & * & 0\vspace*{2pt}\cr
 * & 0 & * & *}\right]$& $S_{2}$\\[20pt]

\includegraphics{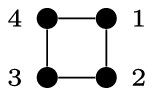}
 & \includegraphics{1130i16} & $\mathbf{T}^{4}$& $D_{4}$\\[3pt]

\includegraphics{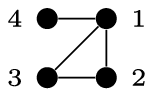}
 & \includegraphics[raise=-12pt]{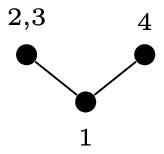} & $\scriptsize\left[\matrix{
 * & 0 & 0 & 0\vspace*{2pt}\cr
 * & * & * & 0\vspace*{2pt}\cr
 * & 0 & 0 & *}\right]$& $\{\mathrm{id}\}$\\[16pt]

\includegraphics{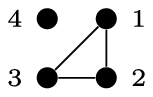}
 &\includegraphics{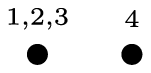} & $\GL_{3}(\RR)\times\mathbf{T}^{1}$& $\{\mathrm{id}
\}$\\[3pt]

\includegraphics{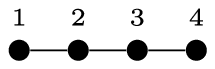}
 &
\includegraphics{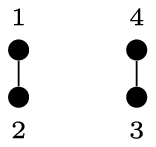}
& $\scriptsize\left[\matrix{
 * & 0 & 0 & 0\vspace*{2pt}\cr
 * & * & 0 & 0\vspace*{2pt}\cr
 0 & 0 & * & *\vspace*{2pt}\cr
 0 & 0 & 0 & *}\right]$& $S_{2}$\\[20pt]

\includegraphics{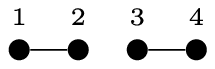}
& \includegraphics{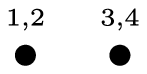} & $\GL_{2}(\RR)\times\GL_{2}(\RR)$& $S_{2}$\\
 \hline
 \end{tabular*}
 \end{table}
\setcounter{table}{0}
 \begin{table}
\caption{(Continued)}
\begin{tabular*}{\textwidth}{@{\extracolsep{\fill}}lccc@{}}
\hline
\multicolumn{1}{@{}l}{$\bolds{\cG}$} & \multicolumn{1}{c}{$\bolds{\bP_{\cC}}$} &
\multicolumn{1}{c}{$\bolds{G^{0}}$} & \multicolumn{1}{c@{}}{$\bolds{\mathrm{Aut}(\widetilde{\cG},c)}$} \\
\hline

\includegraphics{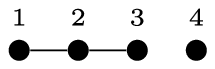}
 &
\includegraphics{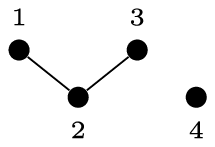}
& $\scriptsize\left[\matrix{
 * & * & 0& 0\vspace*{2pt}\cr
 0 & * & 0& 0\vspace*{2pt}\cr
 0 & * & *& 0\vspace*{2pt}\cr
 0 & 0 & 0 & *}\right]$& $S_{2}$\\[20pt]

\includegraphics{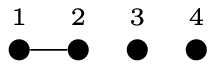}
& \includegraphics{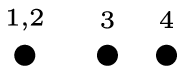} & $\GL_{2}(\RR)\times\mathbf{T}^{2}$&$S_{2}$ \\[3pt]

\includegraphics{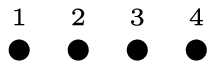}
 & \includegraphics{1130i16}& $\mathbf{T}^{4}$ & $S_{4}$\\
 \hline
 \end{tabular*}
\end{table}


%

\section{What's next} \label{secdiscussion}

In this paper, we presented the complete description of
the maximal subgroup of $GL_m(\RR)$ that stabilizes the Gaussian
graphical model $M(\cG)$ for any given graph $\cG$. The main
motivation for this study was to put Gaussian graphical models into the
framework of (composite) transformation families. Group invariance is a
classical topic in multivariate statistics and there are many ways that
statistical inference can be improved
when the group action is better understood. While we have constructed
the maximal invariant under this group on sample space, we have not yet
used this invariant to develop explicit tests, for example, for model
selection; and while we have given theoretical bounds on when
equivariant estimators for the concentration matrix exist, and how
robust they can be, we have not yet constructed such new estimators. We
regard our work as a step toward achieving these goals for general
graphs, laying down the theoretical framework. On the other hand, in
the case where $G$ acts transitively on the model, we already have a much
better understanding. For instance, it seems feasible to extend the
work of \citet{sun2005estimation} from star-shaped models to general
models in the transitive case. Once these transitive models are
completely understood, it seems natural to move on to those where the
orbit space of $G$ on the model is one-dimensional. Here we expect
beautiful mathematics and statistics to go hand in hand: combinatorics
for characterizing which graphs lead to such models, geometry for a
better understanding of the one-dimensional orbit space, and
statistical inference tailored to the geometry of that space.

\section*{Acknowledgments}
S. Kuhnt thanks the Deutsche Forschungsgemeinschaft (SFB 823, project
B1) for funding.

\begin{supplement}[id=suppA]
\stitle{Proofs and more on the structure of $\mathbf{P}_\cC$}
\slink[doi]{10.1214/13-AOS1130SUPP} 
\sdatatype{.pdf}
\sfilename{aos1130\_supp.pdf}
\sdescription{We provide the proof of Proposition~\ref{propdimstab}
and more results on the structure of
the poset $\mathbf{P}_\cC$ that link our work to \citet{andersson1993}.}
\end{supplement}

%

\printaddresses

\end{document}